\documentclass[11pt]{article}
\usepackage{amssymb,amsfonts,amsmath,latexsym,epsf,tikz,url}

\newtheorem{theorem}{Theorem}[section]

\newtheorem{observation}[theorem]{Observation}
\newtheorem{conjecture}[theorem]{Conjecture}

\newtheorem{lemma}[theorem]{Lemma}

\newcommand{\proof}{\noindent{\bf Proof.\ }}
\newcommand{\qed}{\hfill $\square$\medskip}

\begin{document}

\title{The  distinguishing chromatic number of bipartite graphs of girth at least six}

\author{ Saeid Alikhani  $^{}$\footnote{Corresponding author} 
\and 
Samaneh Soltani
}

\date{\today}

\maketitle

\begin{center}
Department of Mathematics, Yazd University, 89195-741, Yazd, Iran\\
{\tt alikhani@yazd.ac.ir, s.soltani1979@gmail.com}
\end{center}

%%%%%%%%%%%%%%ABSTRACT%%%%%%%%%%%%%%%%%%%%%%%%%%%%%%%%%%%%%%%%%%%%%%%%%%%%%%%%%%%%

\begin{abstract}
The distinguishing number $D(G)$ of a graph $G$ is the least integer $d$
such that $G$ has a vertex labeling   with $d$ labels  that is preserved only by a trivial
automorphism. The distinguishing chromatic number $\chi_{D}(G)$ of $G$ is defined similarly,
where, in addition, $f$ is assumed to be a proper labeling. Motivated by a conjecture in \cite{colins}, we prove that if $G$ is a bipartite graph of girth at least six with the maximum degree $\Delta (G)$,  then    $\chi_{D}(G)\leq \Delta (G)+1$.  We also obtain an upper bound for $\chi_{D}(G)$ where $G$ is a graph with at most one cycle. Finally, we state a relationship between the distinguishing chromatic number of a graph and its spanning subgraphs.

\end{abstract}


\medskip
\noindent{\bf AMS Subj.\ Class.}: 05C25, 05C15

%%%%%%%%%%%%%%%%%%%%%%%%%%%%%%%%%%%%%%%%%%%%%%%%%%%
\section{Introduction }

Let $G=(V,E)$ be a simple finite connected graph. We use the standard graph notation. In particular, ${\rm Aut}(G)$ denotes the automorphism group of $G$.  The \textit{girth} of a graph $G$ is the length of its shortest cycle, and denoted by $g(G)$.  For simple connected graph $G$, and $v \in V$,
the \textit{open neighborhood} of a vertex $v$ is the set $N_G(v) = \{u \in V(G) : uv \in   E(G)\}$. The \textit{closed neighborhood} of a vertex $v$ is the set $N_G[v] = N_G(v)\cup \{v\}$. The \textit{degree} of a vertex $v$ in a graph $G$, denoted by ${\rm deg}_G(v)$, is the number of edges of $G$ incident with $v$. In particular, ${\rm deg}_G(v)$ is the number of neighbours of $v$ in $G$.  We denote by $\Delta(G)$ the  maximum degree of the vertices of $G$. A subgraph $H$ of a graph $G$ is said to be \textit{spanning} if $V(G)=V(H)$ and  $E(H)\subseteq E(G)$. A subgraph $H$ of a graph $G$ is said to be \textit{induced}  if for any pair of vertices $x$ and $y$ of
$H$, $xy$ is an edge of $H$ if and only if $xy$ is an edge of $G$. If $S \subseteq V (G)$ is the vertex set of
$H$, then $H$ can be written as $G[S]$ and is said to be induced by $S$.

 An \textit{$r$-labeling} of the vertices of
a graph $G = (V (G),E(G))$ is a function $f : V (G) \rightarrow \{1, 2, \ldots, r\}$. Then
the labeling $f$ is called a \textit{proper} $r$-labeling if any two adjacent vertices have different
labels. The \textit{chromatic number}  of $G$, denoted by $\chi(G)$, is the minimum $r$ such that $G$ has
a proper $r$-labeling. A concept of symmetry breaking in a graph was introduced by Albertson and
Collins in \cite{Albert}, using the notion of labeling of a graph. An $r$-labeling
is called \textit{distinguishing}  if the trivial automorphism is the only automorphism of $G$
that preserves all the vertex labels. The \textit{distinguishing number}  of $G$, denoted $D(G)$,
is defined to be the minimum $r$ such that $G$ has a distinguishing $r$-labeling. Collins and Trenk \cite{colins}
obtained an analogue of Brook’s theorem. It asserts that  $D(G) \leq \Delta(G) + 1$
holds for any connected graph G. The equality is achieved if and
only if $G = K_{\Delta+1},K_{\Delta,\Delta}$, or $C_5$. They also defined the distinguishing chromatic number which incorporates
the additional requirement that the labeling be proper. They defined the \textit{distinguishing chromatic number}  $\chi_D(G)$ of a graph
$G$ to be the minimum number of labels needed to properly label the vertices of $G$ so
that the only automorphism of $G$ that preserves labels is the identity. Collins and Trenk also proved that $\chi_D(T ) \leq \Delta(T)+1$ for any tree $T$, where the equality is achieved if and only if $T$ belongs to a special
class of trees. In particular, they characterized the connected graphs with maximum
possible distinguishing chromatic number, showing that $\chi_D(G) = |V (G)|$ if and only if $G$
is a complete multipartite graph. Further, they showed that for a connected
graph $G$, $\chi_D(G) \leq 2\Delta$ with equality if and only if $G = K_{\Delta,\Delta}$
or $G = C_6$, a cycle of length six. For connected graphs with $\Delta \leq 2$,
they also completely determined the distinguishing chromatic number.  Laflamme and Seyffarth \cite{Laflamme} stated that if $G$ is
a connected bipartite graph with maximum degree $\Delta\geq 3$, then $\chi_D(G) \leq  2\Delta - 2$
whenever $G \ncong K_{\Delta-1,\Delta}, K_{\Delta,\Delta}$. 

\medskip

In the next section, we improve this bound when $G$ is a connected bipartite graph with girth at least six.
More precisely,  we prove that if $G$ is a bipartite graph of girth at least six with the maximum degree $\Delta (G)$,  then    $\chi_{D}(G)\leq \Delta (G)+1$.  We obtain an upper bound for $\chi_{D}(G)$ where $G$ is a graph with at most one cycle and  a relationship between the distinguishing chromatic number of a graph and its spanning subgraphs in Section 3.

\section{Main result}
 Collins and Trenk \cite{colins} proposed the following conjecture: 
 \begin{conjecture}{\rm\cite{colins}}
 If the girth of a connected graph $G$ is $5$ or greater, then $\chi_D(G) \leq \Delta+1$, where $\Delta \geq 3$.
 \end{conjecture}

  We prove this conjecture for bipartite graphs. For this purpose, we begin with some terminology and background, following \cite{gross}. A \textit{rooted tree}  $(T, z)$ is a tree $T$ with a distinguished vertex $z$, the root. The \textit{depth} or \textit{level}  of a vertex $v$ is its
distance from the root, and the \textit{height}  of a rooted tree is the greatest depth in the tree.
The \textit{parent} of $v$ is the vertex immediately following $v$ on the unique path from $v$ to the root. Two vertices are \textit{siblings}  if they have the same parent.

To prove the main result we state the following observation.
\begin{observation}\label{observ}
Let $G$ be a connected graph with girth at least five.
\begin{enumerate}
\item[(i)]  If  $x,y,z $ are the vertices of $G$ such that $x,y \in N(z)$, then $N[x]\cap N[y]=\{z\}$. 
%In fact, if $w \in N[x]\cap N[y]$, then $wxzyw$ is a 4-cycle, which is a contradiction.

\item[(ii)] Let $x,y,z$ be the vertices of $G$ such that $x,y \in N^i(z)$, where $N^i(z)$ is the set of all vertices of $G$ at distance $i$ of $z$. If $x'\in N^{i+1}(z)\cap N(x)$ and $y_1,y_2\in N^{i+1}(z)\cap N(y)$ such that $y_1 \in N(x')$, then $y_2 \notin N(x')$. 
%In fact, if $y_1,y_2 \in N(x')$, then $x'y_1yy_2x'$ is a 4-cycle, which is a contradiction.

\item[(iii)]  Let $x,y,z,v$ be the vertices of $G$ such that $x,y,z \in N^i(v)$. If $x'\in N^{i+1}(v)\cap N(x)$, $y_1\in N^{i+1}(v)\cap N(y)$ and $z_1\in N^{i+1}(v)\cap N(z)$ such that $y_1 \in N(x')$ and $z_1 \in N(x')$ then $y_1 \notin N(z_1)$. 
%In fact, if $y_1 \in N(z_1)$, then $x'y_1z_1x'$ is a 3-cycle, which is a contradiction.

\item[(iv)]  Let $x,y,z$ be the vertices of $G$ such that $x,y \in N^i(z)$. If $x'\in N^{i+1}(z)\cap N(x)$ and $y'\in N^{i+1}(z)\cap N(y)$, then $|N^{i+1}(z)\cap N(x')\cap N(y')|\leq 1$.  
%In fact, if $w,v\in N(x')\cap N(y')$, then $wx'vy'w$ is a 4-cycle, which is a contradiction.

\item[(v)]  Let $x,y,z$ be the vertices of $G$ such that $x,y \in N^i(z)$. If $x'\in N^{i+1}(z)\cap N(x)$ and $y'\in N^{i+1}(z)\cap N(y)$, and also $x'$ and $y'$ are adjacent, then $ N(x')\cap N(y')=\emptyset$. 
% In fact, if $w\in N(x')\cap N(y')$, then since $x'$ and $y'$ are adjacent so  $wx'y'w$ is a 3-cycle, which is a contradiction.
\end{enumerate}
\end{observation}

\begin{theorem}\label{mainthm}
Let $G$ be a connected bipartite graph with girth at least six, $g(G)\geq 6$. If $\Delta \geq 3$, then  $\chi_D(G) \leq \Delta+1$.
\end{theorem}
\proof  Let $v$ be a vertex of $G$ with maximum degree $\Delta\geq 3$ with $N(v)=\{v_1,...,v_{\Delta}\}$, and $(T,v)$ be a breadth-first search spanning tree rooted at $v$. By a definition of a breadth-first search spanning tree, the distance between $v$ and any vertex $w$ in $T$ is the same as the distance between $v$ and $w$ in $G$. 
%Consider the vertices of $G$ in order they were selected in forming $T$, namely $v,v_1,v_2,\ldots$. Hence $N_G(v)=\{v_1,\ldots , v_{\Delta}\}$ in this ordering. 
 Since $g(G)\geq 6$, so the  induced subgraphs $G[N_G[v_1]\cup N_G[v_2]]$ and $T[N_G[v_1]\cup N_G[v_2]]$ are the same, by Observation \ref{observ}. 
For $1\leq i \leq \Delta$,  we suppose that ${\rm deg}_G(v_i)=t_i+1$, and set  $N_G(v_i)-\{v\}=\{v_{i1},\ldots , v_{it_i}\}$. 
The key point in this proof is that if $N^i(v)$ is the set of all vertices of $G$ at distance $i$ from $v$, then $N_G(x)\cap N^i(v)=\emptyset$ for every $x\in N^i(v)$ and any $1\leq i \leq \Delta$, because $G$ does not have  an odd cycle.
 We state our labeling by the following steps:
\begin{enumerate}
\item[Step 1)] We label the vertex $v$ with 0, and retire the label 0. Thus the vertex $v$ is fixed under each automorphism of $G$ preserving the labeling, and so $N^i(v)$ is mapped to itself, setwise, for every $1\leq i \leq \Delta$, under each automorphism of $G$ preserving the labeling.
 Now we label the vertices $v_1,\ldots , v_{\Delta}$ with labels $1,\ldots , \Delta$, respectively. Then the vertices $v_1,\ldots , v_{\Delta}$  are fixed under each automorphism of $G$ preserving the labeling.
\item[Step 2)]  In this step, we label the vertices $N_G(v_i)-\{v\}$ for every $1\leq i \leq \Delta$, so that each vertex is labeled different from its siblings and its parent in $T$. 
For this purpose, we relate the set $A_i=\{1,2, \ldots , \Delta\}-\{i\}$ of labels to the elements of $N_G(v_i)-\{v\}$ for every $1\leq i \leq \Delta$. Since $g(G)\geq 6$, so for every vertex $x\in N^2(v)$, we have $N(x)\cap N^2(v)=\emptyset$, by Observation \ref{observ}. 
Hence we can label the vertices of  $N_G(v_i)-\{v\}$ with labels in $A_i$, for every $1\leq i \leq \Delta$, so that each vertex is labeled different from its siblings in $T$. 
Now, since none of vertices in $N^2(v)$ are adjacent, so it can be seen that we have  a proper distinguishing labeling of the induced subgraph $G[\{v\}\cup N_G(v)\cup N^2(v)]$ with at most $\Delta +1$ labels $\{0,1,2, \ldots , \Delta\}$.
\end{enumerate}

Before we start to label  the vertices in $N^3(v)$ we state the followin note:
\begin{enumerate}
\item[\bf N1.] For every $1\leq i \leq \Delta$, each permutation of labels  of vertices in  $N_G(v_i)-\{v\}$, or even each new labeling  of these vertices  with labels in $A_i$, so that  each vertex  in  $N_G(v_i)-\{v\}$ is labeled different from its siblings in $T$, makes a new proper distinguishing labeling  of the induced subgraph $G[\{v\}\cup N_G(v)\cup N^2(v)]$ with at most $\Delta +1$ labels.
\end{enumerate}
\begin{enumerate}
\item[Step 3)] Here we label the vertices in $N^3(v)$. For every $1\leq i \leq \Delta$ and $1\leq j \leq t_i$ we define:
\begin{equation*}
M_{ij}:= N_G(v_{ij})\cap N^3(v)=\{w_{(ij)k}~|~1\leq k \leq m_{ij}\}.
\end{equation*}
Hence $|M_{ij}|= m_{ij}$. We want to label the vertices in $M_{ij}$ distinguishingly with labels  $\{1,2, \ldots , \Delta\}$ such that each vertex  is labeled different from its neighbours. For this purpose,   we suppose that for every $1\leq i \leq \Delta$, and $1\leq j \leq t_i$ and  $1\leq k \leq m_{ij}$
$$M_{(ij)k}:=N_G(w_{(ij)k})\cap N^2(v),$$ 
 and $C_{(ij)k}$ is the set of labels are used for elements $M_{(ij)k}$ in Step 2.
  We relate to  vertex $w_{(ij)k}$ the set of labels $B(w_{(ij)k}):=\{1,2,\ldots , \Delta\}\setminus C_{(ij)k}$. 
  If $B(w_{(ij)k})=\emptyset$, then by N1, we can relabel  the vertices of $N^2(v)$ such that $B(w_{(ij)k})\neq \emptyset$. 
 Hence without loss of generality we can assume that $B(w_{(ij)k})\neq \emptyset$, for every  $1\leq i \leq \Delta$, $1\leq j \leq t_i$, and $1\leq k \leq m_{ij}$. Thus the corresponding label set related to the set  $M_{ij}$ is the following set:
  $$B_{ij}=\bigcup_{k=1}^{m_{ij}} B(w_{(ij)k}).$$

 Before we continue our labeling, we state the two following notes:
\begin{enumerate}
\item[\bf N2.] If $|M_{(ij)k}|\geq 2$ for some $1\leq i \leq \Delta$, $1\leq j \leq t_i$, and $1\leq k \leq m_{ij}$, then since $g(G)\geq 6$, so the parents of elements in $M_{(ij)k}$ are different in $T$. 
On the other hand, since $g(G)\geq 6$, so $|N_G(w_{(ij)k})\cap N_G(w')\cap N^2(v)|\leq 1$ for all $w'\in N^3(v)\setminus \{w_{(ij)k}\}$. Hence, the vertex $w_{(ij)k}$ is fixed under each automorphism $f$ of $G$ preserving the labeling, because all vertices in $N^2(v)$ are fixed under $f$. 
Therefore, to have a proper distinguishing labeling, it is sufficient to label the vertices $w_{(ij)k}$ different from the label of its neighbours, i.e., with an arbitrary labels in $ B(w_{(ij)k})$, because $G$ has no odd cycles and hence $N_G(x)\cap N^3(v)=\emptyset$ for every $x\in N^3(v)$.
\item[\bf N3.] Let $M'_{ij}:= \{w_{(ij)k_1}, \ldots , w_{(ij)k_q}\}$ be the set of all elements of $M_{ij}$ for which $|M_{(ij)k}|=1$ where $k\in \{k_1, \ldots ,k_q\}\subseteq \{1,2,\ldots , m_{ij}\}$ for some  $1\leq i \leq \Delta$, $1\leq j \leq t_i$. Then to have a proper distinguishing labeling it is sufficient to label the elements in $M'_{ij}$  different from each other using corresponding label set $B(w_{(ij)k})$ for  $k\in \{k_1, \ldots ,k_q\}$.
\end{enumerate}

Since  $N_G(x)\cap N^3(v)=\emptyset$ for every $x\in N^3(v)$, so by N2 and N3, we can label the vertices of each $M_{ij}$ with corresponding labels in $B_{ij}$, so that the vertices in $M'_{ij}$ have been labeled  different from each other, and also the vertices in $M_{ij}\setminus M'_{ij}$, say $x$, have been labeled with an arbitrary labels in $B(x)$. 
For instance, if $z\in M_{ij}$ and $|N_G(z)\cap N^2(v)|\leq 1$, then we assign the vertex $z$  the smallest label  in $B(z)$ not yet assigned to the already-labeled vertices in $M'_{ij}$. 
If $|N_G(z)\cap N^2(v)|\geq 2$, then we assign the vertex $z$  an arbitrary label in $B(z)$. Therefore, we have a proper distinguishing labeling of the induced subgraph  $G[\{v\}\cup N_G(v)\cup N^2(v)\cup N^3(v)]$, by N2 and N3. 
\end{enumerate}

Before we continue, it must be noted that for every $1\leq i \leq \Delta$, $1\leq j \leq t_i$, each permutation of labels in  $M'_{ij}$, or even each new labeling  of these vertices  with labels in $B_{ij}$, so that  each vertex  in  $M'_{ij}$ is labeled different from each other, makes a new proper distinguishing labeling  of the induced subgraph $G[\{v\}\cup N_G(v)\cup N^2(v)\cup N^3(v)]$. This argument is satisfied for elements in $M_{ij}\setminus M'_{ij}$, too.

For labeling of vertices in $N^i(v)$, $i\geq 4$, we do the similar steps as for $N^3(v)$ and finally, after the finite number of steps, we obtain a proper distinguishing labeling of $G$ with $\Delta +1$ labels. \qed

\section{ Upper bound of $\chi_D(G)$ for trees and unicyclic graphs }
In  this section, we show that if the number of cycles in a connected graph is at most one, then the distinguishing chromatic number is at most $\Delta+1$. 
We need the following lemma:  
\begin{lemma}{\rm \cite{colins}}\label{3.2}
 A labeling of a rooted tree $(T, z)$ in which each
vertex is labeled differently from its siblings and from its parent is a proper distinguishing
labeling.
\end{lemma}

\begin{theorem}
Let $G$ be a connected  graph of order $n$ and size $m$. If $m\leq n$ and $\Delta(G)\geq 3$, then  $\chi_D(G) \leq \Delta+1$.
\end{theorem}
\proof If $m=n-1$, then $G$ is a tree, and we have the result. Then, we suppose that $C$ is the unique cycle in $G$ with alternative vertices $x_1, \ldots , x_t$. It can be seen that $G-E(C)$ is a forest. In fact, $G-E(C)$  is the union of trees $T_{x_i}$, $1\leq i \leq t$, where $T_{x_i}$ has only one common vertex $x_i$ with the cycle $C$.
 Since $\chi_D(C_n)\leq 4$, so we can label the vertices of cycle $C$ with labels 0,1,2, and 3 in a proper distinguishing way. It is clear that the cycle $C$ is mapped to itself under each automorphism of $G$. With respet to the labeling of vertices of $C$, we conclude that the vertices of $C$, i.e., $x_1, \ldots , x_t$, are fixed under each automorphism of $G$ preserving that labeling. To label the vertices of trees  $T_{x_i}$, $1\leq i \leq t$, we note that the vertices  $x_1, \ldots , x_t$ are fixed under each automorphism of $G$ preserving that labeling.
 For every  $T_{x_i}$, we label the adjacent vertices to $x_i$ in  $T_{x_i}$ with ${\rm deg}_G (x_i)-2$ distinct labels, except the label of $x_i$. If these vertices is denoted by $x_{i1}, \ldots , x_{i({\rm deg}_G (x_i)-2)}$, then since the vertices $x_{ij}$, $1\leq j \leq {\rm deg}_G (x_i)-2$,  has at most $\Delta -1$ children, so we can label the children of each of $x_{ij}$, $1\leq j \leq {\rm deg}_G (x_i)-2$,  with distinct labels from the set $\{1,2, \ldots ,\Delta\}$, except the label of vertex $x_{ij}$.  Continuing this process, we can label all vertices of $T_{x_i}$ for $1\leq i \leq t$.  In fact we presented a labeling of the rooted tree $(T_{x_i},x_i)$ for any $1\leq i \leq t$, in which each vertex is labeled different from its siblings and from its parent. Hence, we have a proper distinguishing labeling  of $T_{x_i}$, with $\Delta$ labels, by Lemma \ref{3.2}. Now,  with respect to the   labeling of cycle $C$, we conclude that $\chi_D(G)\leq \Delta+1$. \qed

We end this paper with a theorem about the relationship between the distinguishing chromatic number of  a graph and its spanning subgraph.
\begin{theorem}
Let $G$ be a  graph and $H$ be a spanning subgraph of $G$. If ${\rm Aut}(H)\subseteq {\rm Aut}(G)$, then  $\chi_D(H) \leq \chi_D(G)$.
\end{theorem}
\proof Let $c: V(G)\rightarrow \{1,2, \ldots , \chi_D(G)\}$ be a proper distinguishing labeling of $G$. We claim that the labeling $c$ is a proper distinguishing labeling  of vertices of $H$.  Since $H$ is the spanning subgraph of $G$, so it is clear that $c$ is a proper labeling of $H$. On the other hand, if $f$ is an automorphism of $H$ preserving the labeling $c$, then since  ${\rm Aut}(H)\subseteq {\rm Aut}(G)$, we conclude that $f$ is an automorphism of $G$ preserving the labeling. Since $c$ is a distinguishing labeling of $G$, so $f$ is the identity automorphism. Therefore $c$ is a distinguishing labeling. Now, we have the result.   \qed

\end{document}